\RequirePackage{etoolbox}
\csdef{input@path}{{style/}{graphics/}}
\documentclass[aos,MSNbibl,nameyear,dvips]{arximspdf}

\doi{10.1214/14-AOS1273}% Updated by VTEXPTS2LaTeX.exe, 27.11.2014
\volume{43}
\issue{1}
\pubyear{2015}
\firstpage{299}
\lastpage{322}
\docsubty{FLA}

\makeatletter

\newcommand{\eqref}[1]{(\ref{#1})}
\newcommand{\lVert}{\Vert}
\newcommand{\rVert}{\Vert}
\newcommand{\lvert}{\vert}
\newcommand{\rvert}{\vert}
\newtheorem{theorem}{Theorem}
\newtheorem{lemma}{Lemma}
\newproclaim{example}{Example}
\newproclaim{condition}{Condition}
\newproclaim{remark}{Remark}
\makeatother

\begin{document}
\begin{frontmatter}

\title{Sparsistency and agnostic inference in sparse PCA}
\runtitle{Sparsistency and agnostic inference in sparse PCA}

\begin{aug}
% Corresponding author: Jing Lei - jinglei@andrew.cmu.edu% Updated by
%VTEXPTS2LaTeX.exe, 27.11.2014 13:12
\author[A]{\fnms{Jing}~\snm{Lei}\corref{}\thanksref{t1}\ead
[label=e1]{jinglei@andrew.cmu.edu}\ead
[label=u1,url]{http://www.stat.cmu.edu/\textasciitilde jinglei}}%,
\and
\author[B]{\fnms{Vincent Q.}~\snm{Vu}\ead[label=e2]{vqv@stat.osu.edu}\ead[label=u2,url]{http://vince.vu}}
\runauthor{J. Lei and V. Q. Vu}
\affiliation{Carnegie Mellon University and The Ohio State University}
\address[A]{Department of Statistics\\
Carnegie Mellon University\\
Pittsburgh, Pennsylvania 15213\\
USA\\
\printead{e1}}
\address[B]{Department of Statistics\\
The Ohio State University\\
Columbus, Ohio 43210\\
USA\\
\printead{e2}}
\end{aug}
\thankstext{t1}{Supported by NSF Grant BCS-0941518, NSF Grant
DMS-14-07771 and NIH Grant MH057881.}

% HISTORY:
%
\received{\smonth{1} \syear{2014}}% Updated by VTEXPTS2LaTeX.exe,
%27.11.2014 13:12
%
\revised{\smonth{9} \syear{2014}}% Updated by VTEXPTS2LaTeX.exe,
%27.11.2014 13:12

% ABSTRACT
%
\begin{abstract}
The presence of a sparse ``truth'' has been a constant assumption in the
theoretical analysis of sparse PCA and is often implicit in its
methodological development. This naturally raises questions about the
properties of sparse PCA methods and how they depend on the assumption of
sparsity. Under what conditions can the relevant variables be selected
consistently if the truth is assumed to be sparse? What can be said about
the results of sparse PCA without assuming a sparse and unique truth?
We answer these questions by investigating the properties
of the recently proposed Fantope projection and selection (FPS) method
in the high-dimensional setting. Our results provide general sufficient
conditions for sparsistency of the FPS estimator. These conditions
are weak and can hold in situations where other estimators are known to fail.
On the other hand, without
assuming sparsity or identifiability, we show that FPS provides a sparse,
linear dimension-reducing transformation that is close to the
best possible in terms of maximizing the predictive covariance.
\end{abstract}

% KEYWORDS
\end{frontmatter}

%s1 #&#
\section{Introduction}
\label{sec:introduction}
Sparse principal components analysis (PCA) is a relatively new and popular
technique for simultaneous dimension reduction and variable selection
in high-dimensional data analysis  [e.g., \citet{Jolliffe:2003,Zou:2006}]. It combines
the central idea of classic (or ordinary) PCA
[\citet{Pearson:1901,Hotelling:1933}] with the notion of sparsity: it seeks
linear transformations that reduce the dimension of the data, while depending
on a small number of variables, but retain as much variation as
possible. In
the population setting, these linear transformations correspond to the
projectors of the $k$-dimensional principal subspaces, spanned by the
eigenvectors of the population covariance matrix. The appeal of
sparsity is
that it not only enhances interpretability, but it can yield consistent
estimates when sparsity is truly present in the population, even in high
dimensions [\citet{Johnstone:2009}].

The development of sparse PCA has taken a brisk pace over the past decade.
Methodological developments include regularized estimators based on
penalizing or constraining the variance maximization formulation of PCA
[\citet{Jolliffe:2003,Witten:2009,Journee:2010}],
regression or low-rank approximation [\citet{Zou:2006,Shen:2008}],
convex relaxations [d'Aspremont, Bach and El~Ghaoui (\citeyear{dAspremont:2008jmlr}), \citet{AEJL:2007,fps-nips}],
two-stage procedures based on diagonal thresholding
[\citet{Johnstone:2009,Paul:draft}] and algorithmic variations of
iterative thresholding [\citet{Ma:2013,Yuan:2013}].
Theoretical developments including consistency, rates of convergence,
minimax risk bounds for estimating eigenvectors and principal subspaces
and detection have been established under various statistical models
[\citet{Johnstone:2009,Amini:2009,Lounici:2012,Ma:2013,Berthet:2013a,CMW:2013}, \citeauthor{Vu:2012a} (\citeyear{Vu:2012a,Vu:2013}), \citet{fps-nips}].

The presence of a sparse ``truth'' has been an explicit assumption in the
theoretical analysis of sparse PCA and is often an implicit assumption
in its
methodological development. Here the ``truth'' refers to the leading
$k$-dimensional principal subspace. This naturally raises questions
about the
properties of sparse PCA methods and how they depend on the assumption of
sparsity. Under what conditions can the relevant variables be selected
consistently if the truth is assumed to be sparse? If the truth is not sparse,
and/or not unique, what can be said about the results of sparse PCA?
The first
question is essentially concerned with variable selection consistency, or
\emph{sparsistency}. The second question is a bit more slippery,
because it
essentially requires us to assume nothing beyond independence of the
observations. In other words, the second question is concerned with
\emph{agnostic inference properties} of an estimation method.
In this paper, we investigate variable selection consistency and agnostic
inference properties of the recently proposed \emph{Fantope projection and
selection} (FPS) method due to \citet{fps-nips}.

FPS formulates the sparse PCA problem as a semidefinite program (SDP) whose
solution is a sparse estimate of the projector of the principal
subspace. It
extends the so-called DSPCA formulation of \citet{AEJL:2007} from the
one-dimensional ($k=1$) case to the multidimensional ($k>1$) case, and it
presents a change in perspective by focusing on projectors rather than
individual eigenvectors. FPS is appealing for both theoretical and
computational reasons. Since it directly estimates the projector of the
$k$-dimensional principal subspace, there is no need for iterative
deflation [e.g., \citet{Mackey:2008}], and hence an SDP need only be solved
once rather than $k$ separate times as in DSPCA. \citet{fps-nips}
developed an efficient alternating direction method of multipliers algorithm
to compute FPS, and established $\ell_2$ consistency of FPS under very mild
conditions on the population and input matrices. Most notably, FPS does not
require the stringent spiked covariance model assumption (i.e., the population
covariance matrix is a sparse low-rank matrix plus identity) that is required
by many competing methods such as diagonal thresholding. This makes FPS
applicable to a much wider range of problems, including the important
case of
correlation matrices where diagonal thresholding cannot even be used. (See
Section~\ref{sec:sparsistency} for another example.) However, the variable
selection and agnostic inference properties of FPS remain unknown.

Sparsistency is the ability of an estimator to accurately select the correct
subset of variables when applied to a random sample generated from a model
where only a subset of variables is assumed to be relevant. Conditions under
which sparsistency holds provide important insights about both the
estimator and
the model. They have been studied extensively in other high-dimensional
inference problems such as linear regression
[\citet{Fan:2001,Meinshausen:2006,Zhao:2006,Wainwright:2009}]
and Gaussian
graphical model selection [\citet{Rothman:2008,Lam:2009,Ravikumar:2011}].
In
contrast, theoretical analyses of sparse PCA have mainly focused on consistency
and rates of convergence in matrix norm, with relatively less progress on
variable selection. An exception is \citet{Amini:2009}, who analyzed DSCPCA
under a stringent spiked covariance model with $k= 1$, where the
population covariance matrix is block diagonal and its leading
eigenvector is
assumed to have a small number of nonzero entries of constant
magnitude. Their
work is an important first step, but it leaves open whether or not their
stringent conditions can be loosened and it also does not address the
$k> 1$
case.

In the first part of this paper, we investigate the sparsistency of FPS under
general conditions. Our main results (Theorems \ref{thm:deterministic} and \ref{thm:probabilistic-main})
give broad sufficient conditions under which FPS can
exactly recover the relevant variables. Roughly, the conditions are
that (1)
the relevant variables are not too correlated with the irrelevant variables
(limited correlation), and (2) the leverages (diagonals of the
projector) of
the relevant variables are large enough. Interestingly, these
conditions are
analogous to so-called (1) ``irrepresentability'' and (2) ``$\beta$-min''
conditions for variable selection consistency of the Lasso
[\citet{Zhao:2006,Meinshausen:2006,BuhlmannVanDeGeer}]. To our
knowledge, this is
the first sparsistency result for principal subspaces. When $k= 1$, it
generalizes the results of \citet{Amini:2009} in several directions,
the most
important of which is that it relaxes their block-diagonal condition on the
population covariance matrix.

The second part of this paper addresses the question of assumption-free
interpretation of sparse PCA within a framework that we call
\emph{agnostic inference}. Our goal is to provide both analysis and
interpretation of sparse PCA with essentially no assumptions beyond
independence of observations. The terminology is borrowed from the
learning theory
literature where the chief concern is estimating a classifier or
regression function
without assumptions on the model [\citet{Kearns:1994}];
however, much of
our perspective
is influenced by earlier work on maximum likelihood under misspecification
[\citet{Berk:1966,Huber:1967,White:1982}],
interpretations obtained by
extending the maximum likelihood principle [\citet{Akaike:1973}], and
the notion of \emph{persistence} of high-dimensional linear predictors
proposed by \citet{Greenshtein:2004}. Our point is that although FPS
is derived
under the assumption of sparsity, its results can still be interpreted even
when sparsity does not hold. The main result (Theorem~\ref{thm:agnostic-fps})
is that without assuming sparsity or identifiability, FPS provides a sparse,
linear dimension-reducing transformation that is close to the best possible
in terms of maximizing the \emph{predictive covariance}.

The remainder of the paper is organized as follows. Section~\ref{sec:prelim}
provides the technical background and conditions that are necessary to
state our results---divided between Sections~\ref{sec:sparsistency} (sparsistency)
and \ref{sec:agnostic} (agnostic inference). We discuss these results
in Section~\ref{sec:discussion}, and defer their proofs to  the \hyperref[sec:proof-deterministic]{Appendix}.
Finally, we collect our notation below for our readers' convenience.

\subsection*{Notation}
For two matrices $A,B$ with conformable dimensions, $\langle A, B
\rangle:=\operatorname{trace}(A^T B)$
denotes the trace inner product. For a vector $v\in\mathbb{R}^k$ and
$q\in[0,\infty]$, $\lVert v\rVert_q = {(\sum_{i=1}^p |v_i|^q)}^{1/q}$
is the $\ell_q$ norm if $0<q<\infty$;
when $q=0$, $\lVert v\rVert_0$ is the number of nonzero entries of $v$;
when $q=\infty$, $\lVert v\rVert_\infty=\max_{1\le i\le k} |v_i|$.
For a matrix $A\in\mathbb{R}^{n\times m}$, and index sets
$J_1\subseteq[n]$ and $J_2\subseteq[m]$,
$A_{J_1 J_2}$ denotes the $|J_1|\times|J_2|$ submatrix of $A$
consisting of rows in $J_1$ and columns in
$J_2$, and $A_{J_1 *}$ ($A_{* J_2}$) denotes the submatrix consists of
corresponding rows (columns).
% When $J_1=\{j_1\}$ and/or $J_2=\{j_2\}$, we use simplified notations
%such as $A_{J_1 j_2}$, $A_{j_1 *}$, etc.
Given $q_1, q_2\in[0,\infty]$ and $A\in\mathbb{R}^{n\times m}$, the
matrix $(q_1,q_2)$-pseudonorm $\lVert A\rVert_{q_1,q_2}$ is defined
as ${(\lVert A_{1*}\rVert_{q_1},\lVert A_{2*}\rVert_{q_1},\ldots
,\lVert A_{n*}\rVert_{q_1})}_{q_2}$.
As usual, the spectral norm of $A$ is denoted $\lVert A\rVert$ and the
Frobenius norm is
$\lVert A\rVert_F:=\langle A, A \rangle^{1/2}$.
If $A$ is a symmetric matrix, $\lambda_j(A)$ denotes the $j$th largest
eigenvalue of $A$.
We will use $\Sigma$ to denote the $p\times p$ underlying true
covariance matrix, whose ordered eigenvalues
are $\lambda_1\ge\lambda_2\ge\cdots\ge\lambda_p$.
For a square matrix $A$, $\operatorname{diag}(A)$ denotes its diagonal vector.
For a vector $v$, $\operatorname{supp}(v)$ is the support of $v$ (the
index set corresponding to nonzero entries).

%s2 #&#
\section{Preliminaries}
\label{sec:prelim}%!TEX root = ./agnostic-techreport.tex
% \subsection{The sparse principal subspace problem}

Let $\Sigma\in\mathbb{R}^{p\times p}$ be a symmetric matrix with spectral
decomposition\vspace*{-6pt}
\[
\Sigma= \sum_{j=1}^p
\lambda_j u_j u_j^T,
\]
where $\lambda_1 \geq\cdots\geq\lambda_p$ are eigenvalues and
$u_1,\ldots,u_p \in\mathbb{R}^p$ is an orthonormal basis of
eigenvectors. The $k$-dimensional principal subspace of
$\Sigma$ is the subspace spanned by $u_1,\ldots,u_{k}$. It
is unique if and only if the spectral gap $\lambda_{k} - \lambda
_{k+1} > 0$,
and its projector (orthogonal projection matrix)\vspace*{-3pt} is
\[
\Pi = \sum_{j=1}^{k} u_j
u_j^T = U U^T,
\]
where $U$ is the orthonormal matrix with columns $u_1,\ldots,u_k$.
Every subspace has a unique projector and so we will consider the
principal subspace and $\Pi$ to be equivalent, and
we will also assume that $k$ is known or fixed\vspace*{-3pt} in
advance.

%s2.1 #&#
\subsection{Sparse principal subspaces}\label{sec:sps}
Estimation of the principal subspace requires at minimum that it
be well-defined. When this is the case, we can consider $\Pi$
to be a mapping $x \mapsto\Pi x$ and so it makes sense to
consider indices of the variables that $\Pi$ depends on.
Since $\Pi$ is positive semidefinite, this is equivalent to the
indices of
the nonzero diagonal entries of $\Pi$, because row/column $i$ of
$\Pi$ is nonzero if and only if $\Pi_{ii} \neq0$.

%co1 #&#
\begin{condition}[(SPS)]
$\Sigma$ satisfies the \emph{sparse principal subspace} condition
with support set $J$ if
\renewcommand{\theequation}{SPS}
%
%e1 #&#
\begin{equation}
\label{SPS}
\lambda_{k}(\Sigma)-\lambda_{k+1}(\Sigma) > 0
\quad\mbox{and}\quad \operatorname{supp}\bigl(\operatorname{diag}(\Pi)\bigr) = J.
\end{equation}
\end{condition}

(\ref{SPS}) is the minimal requirement for sparse principal subspace
estimation, and the assumption will only be used in
Section~\ref{sec:sparsistency} in our investigation of sparsistency.
The spectral gap condition ensures that the principal
subspace is identifiable, and the support set definition
states that the principal subspace does not depend on variables outside
of $J$.
This corresponds to a notion of subspace sparsity introduced by \citet{Vu:2013}
called $\ell_0$ \emph{row sparsity}, and it can be shown
that $J = \bigcup_{j=1}^k \operatorname{supp}(u_j)$ for\vspace*{1pt} any orthonormal basis
$\{u_1,\ldots,u_k\}$ of the principal subspace [\citet{Vu:2013}].

%s2.2 #&#
\subsection{Input matrix accuracy}
\label{sec:input-accuracy}
When (\ref{SPS}) is assumed, the main statistical inference problem considered
in this paper is, in a general setting, to estimate $J$ from a
symmetric noisy
version $S$ of $\Sigma$. We then extend the interpretation
and analytical properties of sparse PCA solutions without assuming
(\ref{SPS}).
In both parts,
the estimation accuracy depends on the noisiness of $S$
as an approximation to $\Sigma$, which will be quantified by an
entrywise tail bound on
\[
W:= S- \Sigma.
\]

As motivated by principal component analysis, it may be helpful to
think of~$\Sigma$ as the covariance of a $p$-dimensional
random vector and $S=S_n$ as sample covariance
matrix of a random sample of size $n$, but that is not strictly
necessary for our theoretical analysis. In fact, our sparsistency
results do not even have to assume
that $\Sigma$ or $S$ are positive semidefinite.
In the following, we describe two probabilistic models that imply a
strong entrywise tail bound on $W$.

%ex1 #&#
\begin{example}[(Sample covariance)]
\label{exa:sample-cov}
Let $X, X_1,X_2,\ldots,X_n \in\mathbb{R}^p$ be i.i.d. random
vectors with
$\operatorname{Var}(X) = \Sigma\succeq0$ and let $S$ be the sample
covariance matrix:
\[
S= \frac{1}{n} \sum_{i=1}^n
{(X_i - \bar{X})} {(X_i - \bar{X})}^T,
\]
where $\bar{X} = n^{-1} \sum_{i=1}^n X_i$.
We assume throughout this paper that
\[
\log p \leq n.
\]
By Bernstein's inequality  [\citet{vanderVaartAndWellner},
Chapter~2.2]
if $X$ has sub-Gaussian tails in that there exists constants $K,C>0$
such that
\setcounter{equation}{0}
%e2 #&#
\begin{equation}
\label{eq:X-sub-Gaussian}
\mathbb{P} \bigl(\bigl\lvert v^T(X - \mathbb{E}X)\bigr\rvert\geq t \bigr) \leq K\exp \bigl[-C t^2/\bigl(v^T
\Sigma v\bigr) \bigr]\qquad \mbox{for all }v\neq 0,
\end{equation}
then there is an absolute constant $c > 0$ such that
$S$ satisfies, for $\sigma\ge c \lambda_1$,
%
%e3 #&#
\begin{equation}
\label{eq:max-tail-bound}
\mathbb{P} \biggl( \lVert W\rVert_{\infty,\infty} \geq\sigma\sqrt {
\frac{\log p}{n}} \biggr) \le2p^{-2}.
\end{equation}
In other words, the maximum entrywise error is bounded by $\sigma\sqrt
{\log p/n}$
with high probability. This fact will be the starting point of
subsequent analysis
of the sparsistency of the FPS estimator introduced in Section~\ref{sec:fps}.
The tail bound \eqref{eq:max-tail-bound} is well known and a proof of a
stronger result that implies \eqref{eq:max-tail-bound} can be found in
\citet{Vu:2012a}, Lemma 3.2.2.
\end{example}

%ex2 #&#
\begin{example}[(Random graph models)]\label{exa:clique}
Here, we give an example that does not involve an i.i.d. random sample
and the rate of error bound on $S$ depends only on $p$. Consider a
random graph model with $p$ nodes where edges appear
independently with probability $c_{ij}$ for all $1\le i<j\le p$. Let
$A$ be
the random adjacency matrix such that $A_{ij}=\pm1$ according to the
presence/absence of edge, then the pair $S=AA^T/(p-1)$ and $\Sigma
=\mathbb{E}S$ satisfies
\[
\mathbb{P} \biggl(\lVert S- \Sigma\rVert_{\infty,\infty}\ge c\sqrt {
\frac{\log p}{p-1}} \biggr)\le2p^{-2}
\]
for some universal constant $c$.

This model is related to the planted clique problem
where $c_{ij}=1$ for all $1\le i<j\le s$, and $c_{ij}=1/2$ everywhere else.
The leading eigenvector of
$\Sigma$ is $(1/\sqrt{s},\ldots,1/\sqrt{s},0,\ldots,0)$ and it is supported on
$J=\{1,\ldots,s\}$. Our main result implies that FPS
finds the planted clique with high probability when $s\ge c\sqrt{p\log p}$
for some absolute constant $c$. This is within a factor of $\sqrt{\log p}$
of the best known result for polynomial time recovery in the planted clique
problem [\citet{Deshpande:2013}]. \citet{Berthet:2013b} give another reduction
of the planted clique model to a sparse PCA problem.
\end{example}

For simplicity of notation and presentation, we will focus on the case
of sample covariance matrix in the rest of this paper. But most of our
sparsistency results are applicable to a broader range of problems as
exemplified in the random graph example.

%s2.3 #&#
\subsection{Fantope projection and selection}
\label{sec:fps}
Vu et~al. (\citeyear{fps-nips}) recently proposed an estimator for $\Pi$,
called \emph{Fantope projection and selection} (FPS), defined as
a solution $\widehat{H}$ to the following semidefinite program:
%
%e4 #&#
\begin{equation}
\label{eq:fps}
\widehat{H} :=\arg\max \bigl\{ \langle S, H \rangle- \rho\lVert H
\rVert_{1,1} \bigr\}\quad \mbox{subject to}\quad H\in\mathcal{F}^{k},
\end{equation}
where
\[
\mathcal{F}^{k} := \bigl\{ H\dvtx  0 \preceq H\preceq I
\mbox{ and } \operatorname{trace}(H) = k \bigr\}
\]
is the trace-$k$ \emph{Fantope}, $k > 0$, and $\rho\geq0$
is a tuning parameter. \citet{fps-nips} showed that FPS
can be efficiently computed by alternating direction method of
multipliers
[ADMM, e.g., \citet{Boyd:ADMM}].
When
$\rho=0$, a solution is given by the projector of the $k$-dimensional
principal subspace of $S$ (see Lemma~\ref{lem:fantop-projection} below).
The $\ell_1$ penalty term encourages the solution to be
sparse. Moreover, the decomposability of the $\ell_1$ penalty term
[\citet{NRWY:2012}]
makes it straightforward to analyze the statistical properties of FPS.
In particular, \citet{fps-nips} established a near-optimal Frobenius norm
error bound for the FPS estimator under general conditions.
In the next section, we will show that, if $\Sigma$ satisfies
the (\ref{SPS}) and $S$ satisfies the maximum error bound assumption
\eqref{eq:max-tail-bound}, then under mild conditions,
$\operatorname{supp}[{\operatorname{diag}(\widehat{H})}]=J$ with
high probability for appropriate choices of $\rho$.

In general, the solution to \eqref{eq:fps}, and hence the FPS estimator
may not be unique.
However, we will show that it is unique with high probability when the
(\ref{SPS})
and maximum error bound assumption hold.
The argument utilizes the following elastic net version of FPS:
%
%e5 #&#
\begin{equation}
\label{eq:fps-en}
\max \biggl\{ \langle S, H \rangle - \rho\lVert H
\rVert_{1,1} - \frac{\tau}{2} \lVert H\rVert_{F}^2
\biggr\} \quad\mbox{subject to}\quad H\in\mathcal{F}^{k}.
\end{equation}
Since the objective is a strongly concave function, the solution of
\eqref{eq:fps-en} is unique. A very interesting
and important fact is that when $\rho$ and $\tau$ are small enough,
if a solution of
\eqref{eq:fps} is sparse then it must be the unique solution of \eqref{eq:fps-en}. This observation will
be proved in the  \hyperref[sec:proof-deterministic]{Appendix}  and play a key role in
establishing the uniqueness of
solution for the original FPS problem.

We conclude this section by introducing some basic properties of the
Fantope, which will
be used repeatedly in the proof of main results.
Further properties and discussion of the Fantope will be given in Section~\ref{sec:agnostic}.
Denote the Euclidean projection of a $p\times p$ symmetric matrix $A$ onto
$\mathcal{F}^{k}$ by
\[
\mathcal{P}_{\mathcal{F}^{k}}(A) := \mathop{\arg\min}_{Z\in\mathcal{F}^{k}}\lVert A - Z
\rVert_F^2.
\]

%le1 #&#
\begin{lemma}[(Basic properties of Fantope projection)]\label{lem:fantop-projection}
Let $A$ be a symmetric matrix with eigenvalues $\gamma_1 \geq\cdots
\geq\gamma_p$
and orthonormal eigenvectors $v_1,\ldots,v_p$.
\begin{enumerate}
\item\label{kyfan}
$\max_{H\in\mathcal{F}^{k}} \langle A, H \rangle
= \gamma_1 + \cdots+ \gamma_{k}$ and the maximum is achieved by
the projector of a $k$-dimensional principal subspace of
$A$. Moreover, the maximizer is unique if and only if $\gamma_k>\gamma_{k+1}$.
\item\label{fantopeprojection}
$\mathcal{P}_{\mathcal{F}^{k}}(A) = \sum_j \gamma_j^+(\theta) v_j v_j^T$,
where $\gamma_j^+(\theta) = \min(\max(\gamma_j - \theta, 0), 1)$
and $\theta$
satisfies the equation $\sum_j \gamma_j^+(\theta) = k$.
\item\label{fantopeprojection2}
If $0 < \tau\leq\gamma_k-\gamma_{k+1}$,
then
\[
\mathop{\arg\max}_{H\in\mathcal{F}^{k}}\langle A, H \rangle =
\mathop{\arg\max}_{H\in\mathcal{F}^{k}}\langle A, H \rangle - \frac{\tau}{2} \lVert H
\rVert_F^2 = \mathcal{P}_{\mathcal{F}^{k}}\bigl(
\tau^{-1}A\bigr) = \sum_{j=1}^kv_j
v_j^T,
\]
uniquely.
\end{enumerate}
\end{lemma}

A proof of Lemma~\ref{lem:fantop-projection} is given in Section~\ref{sec:additional-proof}.

%s3 #&#
\section{Sparsistency}
\label{sec:sparsistency}%!TEX root = ./agnostic-techreport.tex
Throughout this section, we assume that $\Sigma$ satisfies
\eqref{SPS} with dimension $k$ and support set $J=\{1,2,\ldots,s\}$
for some $s\ll p$, and that $S$ satisfies the maximum error bound
condition \eqref{eq:max-tail-bound} with some $\sigma>0$. The sample
covariance matrix is covered as a special case in view of Example~\ref{exa:sample-cov}.

Intuitively, variable selection would be easier if
the relevant variables (those in~$J$) and noise variables
(those in $J^c$) are not too
correlated.
In the context of sparse linear regression, such an intuition leads
to the famous Irrepresentable Condition [\citet{Zhao:2006,Wainwright:2009}].
In sparse subspace estimation, we have the analogous Limited Correlation
Condition \eqref{LCC}. In order to state the condition concisely, we
use the following block representation of $\Sigma$:
\[
\Sigma= \pmatrix{
\Sigma_{JJ} &
\Sigma_{J J^c}
\vspace*{3pt}\cr
\Sigma_{J^c J} & \Sigma_{J^c J^c}}.
\]
Similar block representations can be defined for $S$ and $W=S-\Sigma$.

Our main technical condition, the limited correlation condition \eqref{LCC}
is given below.

%co2 #&#
\begin{condition}[(LCC)]
A symmetric matrix $\Sigma$ satisfies the \emph{limited correlation
condition} with constant
$\alpha\in(0,1]$ if
\renewcommand{\theequation}{LCC}
%
%e6 #&#
\begin{equation}
\label{LCC}
\frac{8s}{
\lambda_k(\Sigma)-\lambda_{k+1}(\Sigma)} \lVert\Sigma_{J^c J}
\rVert_{2,\infty} \le 1-\alpha.
\end{equation}
\end{condition}

\eqref{LCC} contains the condition assumed by \citet{Amini:2009} as a
special case, where
$\Sigma_{J^c J}=0$, and hence \eqref{LCC} holds with $\alpha=1$.
Another popular model for sparse PCA is the \emph{spiked covariance
model}, where
$\lambda_k(\Sigma_{JJ})\ge c$, $\Sigma_{J^c J^c}=c I_{p-s}$, and
$\Sigma_{J^c J}=0$.
An important difference between \eqref{LCC} and the assumptions in
previous works is that
previous assumptions, for example, the spiked covariance model, usually
imply that
the relevant variables can be selected with good accuracy by
thresholding the diagonal entries,
while \eqref{LCC} contains situations where such diagonal thresholding
intuition does not work.
Here, we illustrate this difference by a toy example with $p=3$, $k=1$,
$J=\{1,2\}$:
%
%e7 #&#
\setcounter{equation}{4}
\begin{equation}
\label{eq:toy-example}
\Sigma=\pmatrix{
0.9 & 0.8 & t
\cr
0.8 & 0.9 & -t
\cr
t & -t & 1 }.
\end{equation}
This $\Sigma$ satisfies \eqref{LCC} with $\alpha= 0.3$ for any
$|t|\le0.02$, but picking large diagonal entries of $\Sigma$
does not select the relevant variables.

To our knowledge, the \eqref{LCC} is the first sufficient condition for
consistent
sparse PCA variable selection without assuming $\Sigma$ being
block-diagonal and
is also the first sufficient condition for sparse subspace
variable selection consistency.

%s3.1 #&#
\subsection{Sparsistency of FPS}
We state two versions of our main results. The first is a more general,
deterministic
result that provides sufficient conditions for uniqueness, false
positive control,
and false negative control of $\operatorname{supp}(\widehat{H})$. The
second specializes the general
result to the case where $S$ satisfies an entrywise error bound
\eqref{eq:max-tail-bound} like the sample covariance matrix in Example~\ref{exa:sample-cov},
and provides probabilistic guarantees for sparsistency of FPS.

%th1 #&#
\begin{theorem}[(Deterministic support recovery)]\label{thm:deterministic}
Assume $\Sigma$ satisfies \eqref{SPS}.
% Let $\hat Z$ be the solution to FPS problem
% \eqref{eq:fps} with input matrix $\InputMatrix$ and parameters $\rho$.
If the FPS penalty parameter $\rho$ satisfies
%
%e8 #&#
%e9 #&#
\begin{equation}
\label{eq:deterministic-condition-1} \rho^{-1}\lVert S-\Sigma\rVert_{\infty,\infty}+
\frac{8s}{\lambda_k-\lambda_{k+1}}\lVert\Sigma_{J^c J}\rVert _{2,\infty}\le1
\end{equation}
and
\begin{equation}
\label{eq:deterministic-condition-2}
0<\lambda_k-\lambda_{k+1}- 4\rho s \biggl(1+
\frac{8\lambda
_1}{\lambda_k-\lambda_{k+1}} \biggr),
\end{equation}
then the solution $\widehat{H}$ of FPS problem \eqref{eq:fps} is
unique and
satisfies\break $\operatorname{supp}(\operatorname{diag}(\widehat
{H}))\subseteq J$.
If in addition, either
%
%e10 #&#
%e11 #&#
\begin{eqnarray}
\label{eq:signal-strength-1} \min_{j\in J}\sqrt{\Pi_{jj}}& > &
\frac{4\rho s}{\lambda_k-\lambda_{k+1}}\quad \mbox{or }
\\
\label{eq:entry-wise-min}
\min_{(i,j)\in J^2}\lvert\Sigma_{ij}\rvert &>&  2
\rho \quad\mbox{and}\quad \operatorname{rank}\bigl(\operatorname{sign}(
\Sigma_{JJ})\bigr) = 1,
\end{eqnarray}
then the FPS solution satisfies
$\operatorname{supp}(\operatorname{diag}(\widehat{H}))= J$.
\end{theorem}

Theorem~\ref{thm:deterministic} consists of two parts. The first part provides
a set of sufficient conditions [\eqref{eq:deterministic-condition-1} and \eqref{eq:deterministic-condition-2}]
for no false positives. The
second part gives two additional conditions that individually guarantee
no false negatives, and hence exact recovery. We discuss these parts separately.

%pa3.1.subsubsection.1 #&#
\subsubsection*{False positive control}
\eqref{eq:deterministic-condition-1} reveals the motivation for \eqref{LCC}. When
\eqref{LCC} holds, one can choose $\rho=\lVert S-\Sigma\rVert
_{\infty,\infty}/\alpha$
so that \eqref{eq:deterministic-condition-1} holds. On the other hand,
\eqref{eq:deterministic-condition-2} puts some upper bound constraint
on $\rho$.
When $S$ is random and satisfies the maximum error bound condition
\eqref{eq:max-tail-bound}, $\lVert S-\Sigma\rVert_{\infty,\infty}$
depends on
$(n,p,\sigma)$. Then \eqref{eq:deterministic-condition-1} and \eqref{eq:deterministic-condition-2} jointly put
a constraint on $(s,p,n,\sigma,\lambda_1,\lambda_k,\lambda_{k+1})$
so that
there exists a $\rho$ satisfying both conditions.

\eqref{eq:deterministic-condition-2} puts an upper bound on the
sparsity penalty parameter $\rho$. It may seem counterintuitive since
a larger value
of $\rho$ will lead to a sparser solution. In fact, $\rho$
cannot be too large because otherwise the $\ell_1$ penalty term will outweigh
the PCA objective in the FPS problem, leading to a large estimation bias.
Consider the example given in \eqref{eq:toy-example} with $t=0$, if
$S=\Sigma$
and $\rho>0.9$; the FPS solution will return a projection matrix
corresponding to eigenvector $(0,0,1)$, which is supported outside of
the true subset. In general, when $\rho\rightarrow\infty$, the FPS
solution will be a diagonal matrix taking value 1 on diagonal entries
corresponding to the $k$ largest diagonal entries of $S$, and 0 elsewhere.

The proof of false positive control in Theorem~\ref{thm:deterministic}, as
given in Section~\ref{sec:proof-thm:deterministic},
consists of two main steps. The first step (Section~\ref{sec:proof-existence})
is to show that there exists a solution of the FPS problem \eqref{eq:fps} supported on $J$,
using the primal--dual witness (PDW) argument [\citet{Wainwright:2009,Amini:2009,Ravikumar:2011}].
The PDW argument first
constructs a
sparse solution $\tilde H$ supported on $J$ by solving the FPS problem
\eqref{eq:fps}
under additional sparsity constraint $\operatorname
{supp}[\operatorname{diag}(H)]\subseteq J$.
Then it is shown that when $\rho$ is large enough,
with high probability one can find a dual variable $\widehat Z$
such that the primal--dual pair $(\tilde H, \widehat Z)$ satisfies the KKT condition,
and hence is optimal for the original problem. When the solution is
unique, this
ensures that the optimizer is supported on $J$. The challenge here is
to establish KKT
condition when $\Sigma$ is not block diagonal, which requires a
careful and delicate
subspace perturbation analysis in comparing the FPS solution and the
population projector (Lemmas~\ref{lem:sparse-solution} and \ref{lem:primal-dual-optimality}).

The second step is to show that, under the conditions assumed in the
theorem, the
sparse solution constructed in the primal--dual witness argument is
indeed rank-$k$
and also unique. Our proof of uniqueness is novel and makes use of the
elastic net version of FPS \eqref{eq:fps-en}.
% \begin{align}\label{eq:fps-en-1}
% \Estimator(\rho,\tau)=\arg\min_{\PrimalVariable\in\mathcal{F}^{\ndim}} -
% \end{align}
% % The basic properties of Fantope projection (
A key fact used in the proof is that, for
small enough values of $\tau$, the two problems have the same solution
and the uniqueness
of FPS solution follows essentially from that of the elastic net version.
The details are given in Section~\ref{sec:uniqueness-proof}. % The conclusion
%that
% $\supp[\diag(\tilde Z)]=J$ can be reached by a Frobenius norm error
%bound argument.

%pa3.1.subsubsection.2 #&#
\subsubsection*{False negative control}
Having established false positive control in Theorem~\ref{thm:deterministic},
full sparsistency will be established if we can show that the number of
false negatives is also zero.
In sparsity pattern recovery, the number of false negatives is
typically controlled by assuming
a lower bound on the magnitude of signals carried by relevant
variables. In the context of
principal subspace estimation, our first sufficient condition for false
negative control \eqref{eq:signal-strength-1} originates from a
Frobenius norm error bound of FPS established in \citet{fps-nips}:
%
%e12 #&#
\begin{equation}
\label{eq:frobenius-error}
\lVert\widehat{H}-\Pi\rVert_F\le\frac{4\rho s}{\lambda_k-\lambda
_{k+1}} .
\end{equation}

The other sufficient condition for controlling false negative \eqref{eq:entry-wise-min} is
motivated by an assumption used by \citet{Amini:2009} for the $k=1$ case
where the leading eigenvector is assumed to be $v_1=(\mathbf
{1}_s/\sqrt{s},0)$
(where $\mathbf{1}_s$ is the $s\times1$ vector of ones, and the signs
of nonzero entries can actually be arbitrary)
and $\Sigma_{JJ}=\theta v_1v_1^T+ I_{s}$.
Let $\operatorname{sign}(\Sigma_{JJ})$ be the $s\times s$ matrix of
entry-wise signs of $\Sigma_{JJ}$.
%using a direct characterization
%of the leading eigenvector of $\PopMatrix_{JJ}$.
% A proof of \Cref{thm:entry-wise} is given in
Our condition \eqref{eq:entry-wise-min} generalizes that of \citet{Amini:2009} in three directions. First,
we allow principal subspaces of dimension $k> 1$. Second, we allow nonzero
correlation between the relevant and irrelevant variables, whereas
\citet{Amini:2009} assumes a block diagonal structure. Third, we do not
require a generalized spiked covariance model as in \citet{Amini:2009}.
The proof of the second part of Theorem~\ref{thm:deterministic} is given in Section~\ref{sec:additional-proof}.

%th2 #&#
\begin{theorem}[(Sparsistency)]
\label{thm:probabilistic-main}
Assume that $\Sigma$ satisfies \eqref{SPS} and \eqref{LCC},
and that $S$ satisfies
the maximum error bound \eqref{eq:max-tail-bound} with scaling factor
$\sigma$.
If
%
%e13 #&#
\begin{equation}
\label{eq:sample-complexity}
s\sqrt{\frac{\log p}{n}}< \frac{\alpha{(\lambda_k-\lambda_{k+1})}^2}{
4\sigma(8\lambda_1+\lambda_k-\lambda_{k+1})},
\end{equation}
and the FPS penalty parameter $\rho$ in \eqref{eq:fps} satisfies
\[
\rho= \frac{\sigma}{\alpha}\sqrt{\frac{\log p}{n}},
\]
then with probability at least $1-2p^{-2}$, the FPS estimate $\widehat
{H}$ is
unique and satisfies $\operatorname{supp}(\operatorname
{diag}(\widehat{H}))\subseteq J$.
If in addition, either
%
%e14 #&#
%e15 #&#
\begin{eqnarray}
\label{eq:lower-bound-condition-prob-1}
\min_{j\in J}\sqrt{\Pi_{jj}} &>&
\frac{4 s \sigma}{\alpha(\lambda_k-\lambda_{k+1})} \sqrt{\frac{\log p}{n}} \quad\mbox{or}
\\
\label{eq:lower-bound-condition-prob-2}
\min_{(i,j)\in J^2}\Sigma_{ij} &>&
\frac{2\sigma}{\alpha} \sqrt{\frac{\log p}{n}}\quad \mbox{and}\quad \operatorname{rank}\bigl(
\operatorname{sign}(\Sigma_{JJ})\bigr) = 1,
\end{eqnarray}
then $\operatorname{supp}(\operatorname{diag}(\widehat{H})) = J$.
\end{theorem}

\begin{pf}
%{Proof of Theorem~\ref{thm:probabilistic-main}}
Using the maximum error bound condition, with probability at least
$1-2p^{-2}$ we have
$\rho^{-1}\lVert S-\Sigma\rVert_{\infty,\infty}\le\alpha$. This
together with
the property \eqref{LCC} of~$\Sigma$ establishes \eqref{eq:deterministic-condition-1}.
On the other hand, \eqref{eq:sample-complexity}
ensures that \eqref{eq:deterministic-condition-2} holds.
On the other hand, \eqref{eq:lower-bound-condition-prob-1} implies
\eqref{eq:signal-strength-1}, and \eqref{eq:lower-bound-condition-prob-2} implies
that
the choice of $\rho$
satisfies \eqref{eq:entry-wise-min}. The claimed results follow from
Theorem~\ref{thm:deterministic}.
\end{pf}

%re1 #&#
\begin{remark}
When the eigenvalues of $\Sigma$ are constants and do
not change with $(n,p,s)$,
Theorem~\ref{thm:probabilistic-main}
recovers a rate developed by \citet{Amini:2009} as a special case where
%$d=1$
%and $\InputMatrix$ is the sample covariance with a block diagonal $
%leading eigenvector is $(\onevec_s,0)/\sqrt{s}$ (the signs can be
%arbitrary). In this case,
Theorem~\ref{thm:probabilistic-main} implies that a sufficient condition
for consistent variable selection (with suitable choice of $\rho$)
is $s\sqrt{\log p/n}\le c$ for a constant $c$ [according to \eqref{eq:sample-complexity}
and \eqref{eq:lower-bound-condition-prob-2}]. \citet{Amini:2009} also obtain
a sharper sufficient condition $s\log p/n\le c'$, by assuming that the
solution is rank 1. However,
\citet{Krauthgamer:2013} show that, with high probability, the solution
is not rank 1 unless
$s\sqrt{1/n}$ is bounded by a constant.
\end{remark}

%re2 #&#
\begin{remark}
Condition \eqref{eq:sample-complexity} suggests that the required
sample size needs to increase as $\lambda_1$ increases. This is because
the oracle operator norm error bound of the principal subspace (i.e.,
assuming $J$ is known) has a factor of $\lambda_1$. In an extremal
case, when $\lambda_1$ is large and $\lambda_j$ ($j\ge2$) are much
smaller, the estimation error of the leading eigenvector will likely
dominate all the remaining spectral gaps, making it hard to recover the
remaining eigenvectors.
\end{remark}

%s4 #&#
\section{Agnostic inference}
\label{sec:agnostic}%!TEX root = ./agnostic-techreport.tex
$\!$Consistent estimation and variable selection inevitably depend
on the existence of a ``true'' model. For sparse PCA,
this corresponds to the assumption that the $k$-dimensional
principal subspace of $\Sigma$ is (1) identifiable
and~(2) sparse. Under this assumption, previous work [e.g.,
\citet{fps-nips}]
and the theory presented in Section~\ref{sec:sparsistency} establish
conditions under
which consistent estimation and variable selection are possible.
While these results can provide useful insights for sparse PCA and
FPS, the conditions may or may not hold in practice. Therefore,
it is important to understand the statistical inference problem
without these assumptions. This is the \emph{agnostic inference}
perspective. Can we remove the assumptions of identifiability
and sparsity? Is there an assumption-free interpretation for FPS?

Without assuming identifiability, variable selection and estimation
consistency are no longer valid objectives, since there is no unique
``true'' parameter to estimate. For example, when $\Sigma=I$,
every $k$-dimensional subspace is a principal subspace, and even
if there is a unique principal subspace, it may not be sparse.
To develop an assumption-free interpretation, we return to the
basic objective function of PCA. Let $X$ be a random vector with
covariance matrix $\Sigma$. PCA can be interpreted as a covariance
maximization technique. It seeks a rank-$k$
projector $H$ that maximizes the
\emph{predictive covariance}:
\[
% \label{eq:fps-objective}
\operatorname{trace}\bigl( \operatorname{Cov}(X, HX | H) \bigr) =
\langle\Sigma, H \rangle.
\]
If we interpret $H$ as a dimension-reducing transformation,
then $\langle\Sigma, H \rangle$ is just the total covariance
between the input $X$ and output $HX$.

%s4.1 #&#
\subsection{Sparse and shrinking dimension reduction}
FPS also maximizes covariance, but it replaces the rank-$k$ projector
constraint on $H$ with a Fantope constraint and
an additional sparsity constraint via the $(1,1)$-norm. Let
%
%e16 #&#
\begin{equation}
\label{eq:fps-constraint} \widehat{H}_R := \mathop{\arg\max}_{H\in\mathcal{F}^{k}, \lVert H\rVert_{1,1} \leq R}
\langle S, H \rangle.
\end{equation}
By Lagrangian duality, this constrained form of FPS is equivalent to
the penalized form \eqref{eq:fps} in the sense that given $S$,
for every $R$ there is a corresponding $\rho$ such that a solution
of \eqref{eq:fps} is also a solution of \eqref{eq:fps-constraint} and
vice-versa. The corresponding population version of \eqref{eq:fps-constraint}
is
%
%e17 #&#
\begin{equation}
\label{eq:fps-constraint-pop}
H_R := \mathop{\arg\max}_{H\in\mathcal{F}^{k},\lVert H\rVert_{1,1}\leq R} \langle
\Sigma, H \rangle.
\end{equation}

The meaning of $H\in\mathcal{F}^{k}$ may be unclear since it
is not necessarily a rank-$k$ projector. However, it turns out that if
we regard $H$ as a linear transformation
$x \mapsto Hx$, then $H$ is a
\emph{smoother matrix} [\citet{ESL}, Section~5.4.1]  and the Fantope
coincides with
a class of linear smoothers called \emph{shrinking smoothers}
[\citet{Buja:1989}].
The two essential properties of $H$ are:
\begin{enumerate}
\item
$0 \preceq H\preceq I$. This is equivalent to the
condition that
\[
\lVert x\rVert^2 \geq\lVert Hx\rVert^2 + \lVert x - Hx
\rVert^2 \qquad\mbox{for all $x$.}
\]
In other words, the sum of squares of the transformation $Hx$
and its residual $x - Hx$ cannot be larger than that of $x$.
A map satisfying this property is called \emph{firmly nonexpansive}.
\item$\operatorname{trace}(H) = k$. If $H$ is a projector,
then $k$ is the dimension of the projection space. It is also equal to
$\operatorname{trace}[\operatorname{Cov}(\xi, H\xi)]$ when $\xi$
is a random vector with
$\operatorname{Var}(\xi) = I$. By analogy,
$\operatorname{trace}(H)$ is the \emph{effective degrees of freedom} of
$H$ [see \citet{ESL},  Section~5.4.1].
\end{enumerate}
These two properties are exactly those laid out by \citet{ESL} for
smoother matrices and shrinking smoothers. In the context of dimension
reduction,
we call the action of $H\in\mathcal{F}^{k}$
\emph{shrinking dimension reduction}.

Now we turn to the $(1,1)$ norm constraint in \eqref{eq:fps-constraint-pop}.
A natural notion of sparsity of a matrix $H\in\mathcal{F}^{k}$
is $\lVert H\rVert_{2,0}$, the number
of nonzero rows. Here, we use the $(1,1)$-norm as an alternative convex
measure of sparsity. For $H\in\mathcal{F}^{k}$ we have,
by Cauchy--Schwarz,
%
%e18 #&#
\begin{equation}
\label{eq:1-1-bounded-by-2-0}
\lVert H\rVert_{1,1}\le k \lVert H\rVert_{2,0}.
\end{equation}
That is, if $\lVert H\rVert_{2,0}$ is small, then
$\lVert H\rVert_{1,1}$ must also be small.
% See FPS long version for further discussion of the $(1,1)$-norm bound
%and
% the Fantope

%s4.2 #&#
\subsection{Persistence of FPS}
Our main result in the assumption-free setting is an interpretation of the
constrained form of FPS and its persistence under no assumptions on
$\Sigma$.

%th3 #&#
\begin{theorem}[(Persistence)]
\label{thm:agnostic-fps}
Let $X,X_1,\ldots,X_n \in\mathbb{R}^p$ be i.i.d. random vectors that
satisfy the tail probability bound \eqref{eq:X-sub-Gaussian} (i.e.,
$X$ is sub-Gaussian). Then with probability at least $1-2p^{-2}$,
\[
\langle\Sigma, H_R \rangle \geq \langle\Sigma,
\widehat{H}_R \rangle \geq \langle\Sigma, H_R \rangle -
c R \lambda_1 \sqrt{\frac{\log p}{n}},
\]
where $c > 0$ is a constant.
\end{theorem}

Our proof of Theorem~\ref{thm:agnostic-fps} is given in Section~\ref{sec:additional-proof}.
Theorem~\ref{thm:agnostic-fps} shows that the predictive covariance of FPS comes
close to that of the best sparse $H$ in the Fantope.
This is essentially an assumption-free interpretation.

Let
\[
\Pi_{k,s} := \arg\max \bigl\{ \langle\Sigma, \Pi\rangle\dvtx \Pi\mbox{ is a
rank-$k$ projector and $\lVert\Pi\rVert_{2,0} \leq s$}\bigr\}
\]
be the best rank-$k$ and $s$-sparse projector.
What can we say about $\widehat{H}_R$ and $\Pi_{k,s}$?
In this case, \eqref{eq:1-1-bounded-by-2-0} implies that $\lVert\Pi
_{k,s}\rVert_{1,1}\le ks$. Thus
$\Pi_{k,s}$ is in the feasible set of \eqref{eq:fps-constraint}
if $R\ge ks$.
If we do not assume any structure on $\Sigma$,
Theorem~\ref{thm:agnostic-fps} implies that, with high probability,
\[
\langle\Sigma, \widehat{H}_R \rangle \geq \langle\Sigma,
\Pi_{k,s} \rangle - cR\lambda_1\sqrt{\frac{\log p}{n}},
\]
when $R\ge ks$. If we assume in addition that $\Sigma$
does have a $k$-dimensional principal subspace involving at most $s$ variables,
then the result can be strengthened to
\[
\langle\Sigma, \Pi_{k,s} \rangle \geq \langle\Sigma,
\widehat{H}_{R} \rangle \geq \langle\Sigma, \Pi_{k,s} \rangle
- cR\lambda_1\sqrt{\frac{\log p}{n}}.
\]
Here, the assumption that $\Sigma$ has a sparse principal subspace is still
much weaker than the sparse principal subspace condition required by the
sparsistency argument in Section~\ref{sec:sparsistency}, because there is no
requirement on uniqueness of the principal subspace. As a simple example,
$\Sigma= I$ satisfies the sparsity condition but
not the uniqueness condition.

%re3 #&#
\begin{remark}[(Stability of FPS)]
A referee has pointed out to us that there is another interpretation of
Theorem~\ref{thm:agnostic-fps} in terms of the \emph{continuity} of the maximal
predictive covariance map
\[
f(\Omega) :=\max\bigl\{ \langle\Omega, H \rangle\dvtx  H\in\mathcal{F}^{k},
\lVert H\rVert_{1,1} \leq R \bigr\}.
\]
The proof of Theorem~\ref{thm:agnostic-fps} implies that
\[
\bigl\lvert f(\Sigma+ \Delta) - f(\Sigma) \bigr\rvert \leq 2 R \lVert\Delta
\rVert_{\infty,\infty}.
\]
So the predictive covariance of FPS is relatively stable under perturbations
of $\Sigma$ if $R \lVert\Delta\rVert_{\infty,\infty}$ is small.
\end{remark}

%s5 #&#
\section{Discussion}
\label{sec:discussion}%!TEX root = ./agnostic-techreport.tex
A connection between sparse PCA and sparse linear regression has been observed
by \citet{Vu:2013}. They established minimax rates for estimation under
$\ell_2$ loss with $\ell_q$-penalized estimators with suitably
defined model
parameters and observed that the rates are identical to those for sparse
linear regression when the effective noise variance is defined appropriately.
The sparsistency result in the present paper further extends this connection
to variable selection. Roughly speaking, the previously used spiked
covariance model in sparse PCA, which assumes that
\[
\Sigma= U\Lambda U^T + \sigma^2 I_{p},
\]
where $U$ is $p \times k$ orthonormal matrix and $\Lambda\succ0$ is
diagonal  [see, e.g., \citet{Johnstone:2009,Birnbaum:2012,Ma:2013,CMW:2013}],
corresponds to the orthogonal design in linear regression, in the sense that
the relevant and noise variables are not correlated. Moreover, the
$\sigma^2I$ term boosts the signal by adding $\sigma^2$ to all the relevant
diagonal entries in $\Sigma$ and, therefore, thresholding based methods
usually work well. The limited correlation condition developed in this paper
is analogous to the irrepresentable condition
[\citet{Zhao:2006,Meinshausen:2006}] for $\ell_1$-penalized sparse regression
(Lasso), where convex optimization methods can succeed when the correlation
between relevant and noise variables is small.

When the eigenvalues of $\Sigma$ are fixed, a sufficient condition for
consistent variable selection using FPS is $s\lesssim\sqrt{n/\log
p}$. This
is comparable to the corresponding rate developed for $k=1$ by
\citet{Amini:2009} when the rank of the solution is not assumed to be $1$.
It has been shown by \citet{Amini:2009} that the information-theoretic
critical rate is $s\lesssim n/\log p$. That is, if $s\gg n/\log p$, no method
can succeed in variable selection. It remains an open question if there exist
polynomial time methods that can consistently select relevant variables
in the
range $\sqrt{n/\log p}\lesssim s\lesssim n/\log p$. An interesting
work in
this direction is that by \citet{Berthet:2013b}, which shows that, for
$k=1$, testing a sparse PCA model in this regime is at least as hard as
solving the planted clique problem beyond the well-believed computational
barrier.

The predictive covariance maximization interpretation of PCA leads to a
natural characterization of the Fantope as the collection of all shrinking
smoothers with $k$ effective degrees of freedom.
Without any assumptions on $\Sigma$, FPS gives us a dimension reducing
transformation that is sparse while being computationally tractable,
and it
nearly approaches the best predictive covariance. In practice, it would be
useful to estimate the predictive covariance of the FPS solution for a
particular value of $\rho$ using risk estimates such as
cross-validation. This leads to a data-driven procedure for selecting
the best FPS
tuning parameter $\rho$. The detailed design and properties of such
a cross-validation method is an important and interesting topic for future
work.

\begin{appendix}\label{sec:proof-deterministic}
%s6 #&#
\section*{Appendix: Technical proofs}
%!TEX root = ./agnostic-techreport.tex

This appendix  contains detailed technical proofs. In
Section~\ref{sec:proof-thm:deterministic}, we prove the deterministic sparsistency
theorem (Theorem~\ref{thm:deterministic}). Other proofs, including those of
Lemma~\ref{lem:fantop-projection} and
Theorem~\ref{thm:agnostic-fps} are given in Section~\ref{sec:additional-proof}.

%s6.1 #&#
\subsection{Proof of Theorem~\texorpdfstring{\protect\ref{thm:deterministic}}{1}}
\label{sec:proof-thm:deterministic}
%s6.1.1 #&#
\subsubsection{Existence of a sparse solution}
\label{sec:proof-existence}
The primal--dual witness argument starts from the dual form of the FPS
problem \eqref{eq:fps}.
Using strong duality, we can write \eqref{eq:fps} in a equivalent
min--max form:
%
%e19 #&#
%e20 #&#
\begin{eqnarray}
&& \max_{H\in\mathcal{F}^{k}} \langle S, H \rangle- \rho\lVert H
\rVert_{1,1}
\nonumber
\\[1.5pt]
&& \qquad \iff\quad\max_{H\in\mathcal{F}^{k}}\min_{Z\in\mathbb{B}_p} \langle S, H
\rangle-\rho\langle H, Z \rangle-k\rho
\nonumber
\\[1.5pt]
\label{eq:fps-min-max} && \qquad \iff\quad\max_{H\in\mathcal{F}^{k}}\min_{Z\in\mathbb{B}_p}
\langle S- \rho Z, H \rangle
\\[1.5pt]
\label{eq:fps-max-min}
&& \qquad\iff\quad\min_{Z\in\mathbb{B}_p}\max_{H\in\mathcal{F}^{k}} \langle S- \rho
Z, H \rangle,
\end{eqnarray}
where\vspace*{1pt}
$\mathbb{B}_p=\{Z\in\mathbb{R}^{p\times p}\dvtx \operatorname
{diag}(Z)=0,Z=Z^T,\lVert Z\rVert_{\infty,\infty}\le1\}$.
According to the standard Karush--Kuhn--Tucker (KKT) condition, a pair
$(\widehat{H},\widehat{Z})\in\mathcal{F}^{k}\times
\mathbb{B}_p$ is optimal for problems \eqref{eq:fps-min-max} and
\eqref{eq:fps-max-min} if and only if
%
%e21 #&#
%e22 #&#
%e23 #&#
\begin{eqnarray}\label{eq:kkt-1}
\widehat{Z}_{ij} &=&  \operatorname{sign}(\widehat{H}_{ij})
\qquad \forall i\neq j, \widehat{H}_{ij}\neq0,
\\[1.5pt]
\label{eq:kkt-2}
\widehat{Z}_{ij} &\in& [-1,1]\qquad \forall i\neq j, \widehat
{H}_{ij}=0,
\\[1.5pt]
\label{eq:kkt-3}
\widehat{H}&=& \mathop{\arg\max}_{H\in\mathcal{F}^{k}} \langle S-\rho \widehat{Z}, H \rangle.
\end{eqnarray}

To proceed with the primal--dual witness argument, we first construct an
additionally constrained solution $\tilde H$
as follows:
%
%e24 #&#
\begin{equation}
\label{eq:sparse-solution}
\tilde H=\mathop{\arg\max}_{H\in\mathcal{F}^{k},\operatorname
{supp}(\operatorname{diag}(H))\subseteq J} \langle S, H \rangle-
\rho\lVert H\rVert_{1,1}.
\end{equation}
Let $\tilde Z$ be a corresponding optimal dual variable. By Lemma~\ref{lem:sparse-solution},
$\tilde H$ is a rank-$k$ projector supported on $J$.

Let ${\hat{U}_J\choose 0}$ and ${U_J\choose 0}$ be $p\times k$ orthogonal matrices consisting of the $k$ leading
eigenvectors of $S-\rho\tilde Z$ and $\Sigma$, respectively,
where $\hat U_J$ and $U_J$ are $s\times k$ orthogonal matrices.
According to Lemma~\ref{lem:sparse-solution},
there exists
a $s\times s$ orthonormal matrix $Q$ such that\vspace*{1pt}
$\hat U_J = Q U_J$ and $\lVert Q-I\rVert_F\le8\rho s/(\lambda
_k-\lambda_{k+1})$.

Define a modified primal--dual pair $(\widehat{H},\widehat{Z})$ as
follows (recall that $W= S-\Sigma$):
%
%e25 #&#
%e26 #&#
%e27 #&#
\begin{eqnarray}
\widehat{H}& = &\tilde H,
\nonumber
\\[1.5pt]
\label{eq:dual-construct-1}
\widehat{Z}_{JJ}&=&\tilde Z_{JJ},
\\[1.5pt]
\label{eq:dual-construct-2}
\widehat{Z}_{ij}&=&\frac{1}{\rho} \bigl\{S_{ij}-\langle
Q_{i*}, \Sigma_{J,j} \rangle \bigr\}, \qquad  (i,j)\in J\times
J^c,
\\[1.5pt]
\label{eq:dual-construct-3}
\widehat{Z}_{ij} &=&\frac{1}{\rho}W_{ij}, \qquad (i,j)\in
\bigl(J^c\bigr)^2, i\neq j.\\[-30pt]\nonumber
\end{eqnarray}\vadjust{\goodbreak}

We need to check that $(\widehat{H}, \widehat{Z})$ is feasible for
(\ref{eq:fps-min-max}) and (\ref{eq:fps-max-min}) and satisfies the KKT conditions
(\ref{eq:kkt-1}) to (\ref{eq:kkt-3}).

\textit{Checking feasibility}. The feasibility of $\widehat{H}$ is
obvious. To
check feasibility of $\widehat{Z}$, we only need to verify
that $\widehat{Z}_{ij}\in[-1,1]$ for all $(i,j)\in J\times J^c$. In fact,
\begin{eqnarray*}
|\widehat{Z}_{ij}| & \le & \frac{1}{\rho} \bigl[ |S_{ij}-
\Sigma_{ij}|+\bigl|\Sigma_{ij}-\langle Q_{i*},
\Sigma_{J,j} \rangle\bigr| \bigr]
\\[1.5pt]
& \le& \frac{1}{\rho} \bigl[\lVert W\rVert_{\infty,\infty} +\bigl
\lVert(I-Q)_{i*}\bigr\rVert\times\lVert\Sigma_{J,j}\rVert
\bigr]
\\[1.5pt]
& \le& \frac{1}{\rho} \bigl[\lVert W\rVert_{\infty,\infty} +\lVert I-Q
\rVert_F \lVert\Sigma_{J^c J}\rVert_{2,\infty} \bigr]
\\[1.5pt]
& \le& \frac{1}{\rho} \biggl[\lVert W\rVert_{\infty,\infty}+ \frac{8\rho s}{\lambda_k-\lambda_{k+1}}
\lVert\Sigma_{J^c J}\rVert _{2,\infty} \biggr] \le 1,
\end{eqnarray*}
where the last inequality follows from
\eqref{eq:deterministic-condition-1}.

\textit{Checking KKT condition} \eqref{eq:kkt-1}. Because $\widehat
{H}$ only has nonzero entries in
$J\times J$, so $(\widehat{H},\widehat{Z})$
satisfies \eqref{eq:kkt-1} by construction.

\textit{Checking KKT condition} \eqref{eq:kkt-2}. For $(i,j)$ in
$J\times J$, \eqref{eq:kkt-2} is satisfied
for $(\widehat{H},\widehat{Z})$ because the same condition is
satisfied for $(\tilde H,\tilde Z)$. For
$(i,j)\notin J\times J$, we have $\widehat{H}_{ij}=0$ and \eqref{eq:kkt-2} follows from the feasibility
of $\widehat{Z}$.

\textit{Checking KKT condition} \eqref{eq:kkt-3}.
Recall that $W=S-\Sigma$. Let $\tilde W$
be the $(p-s)\times(p-s)$ diagonal matrix that agrees with $W_{J^c J^c}$
on diagonal entries.
By Lemma~\ref{lem:fantop-projection}, it suffices to show
that ${\hat{U}_J\choose 0}$ spans a $k$-dimensional principal subspace of
%
%e28 #&#
\begin{equation}
\label{eq:sigma-tilde}
\tilde\Sigma:= S-\rho\widehat{Z}= \pmatrix{ %
S_{JJ}-
\rho\tilde Z_{JJ} & Q \Sigma_{JJ^c}
\vspace*{3pt}\cr
\Sigma_{J^c J}Q^T
& \tilde W+\Sigma_{J^c J^c}},
\end{equation}
which is established in Lemma~\ref{lem:primal-dual-optimality}.

Now we have shown that $(\widehat{H},\widehat{Z})$ is indeed an optimal
primal--dual pair for \eqref{eq:fps-min-max} and \eqref{eq:fps-max-min}, and hence
$\widehat{H}$ is a solution of \eqref{eq:fps} and is also supported
only on $J$.

%s6.1.2 #&#
\subsubsection{Uniqueness of solution}\label{sec:uniqueness-proof}
Consider the elastic net version of FPS in~\eqref{eq:fps-en} and its
max--min and min--max forms
using dual variable $Z\in\mathbb{B}_p:=\{Z\in\mathbb
{R}^{p\times p}\dvtx \operatorname{diag}(Z)=0,Z=Z^T,\lVert Z\rVert_{\infty
,\infty}\le1\}$:
\begin{eqnarray*}
\label{eq:fps-min-max-en}
&& \min_{H\in\mathcal{F}^{d}}\max_{Z\in\mathbb{B}_p} -
\langle S, H \rangle+\rho\langle H, Z \rangle+\frac{\tau
}{2}\lVert H
\rVert_F^2
\\
&& \qquad \iff \quad \max_{Z\in\mathbb{B}_p}\min_{H\in\mathcal{F}^{k}}
\frac{\tau}{2}\biggl\lVert  {H- \frac{1}{\tau} (S- \rho Z
)}\biggr\rVert_F^2 - \frac{1}{2\tau}\lVert S- \rho Z
\rVert_F^2.
\end{eqnarray*}
\noindent The KKT condition for optimality of $(\widehat{H},\widehat{Z})\in
\mathcal{F}^{k}\times
\mathbb{B}_p$ becomes
%
%e29 #&#
%e30 #&#
%e31 #&#
\begin{eqnarray}
\widehat{Z}_{ij} &=&  \operatorname{sign}(\widehat{H}_{ij})
\qquad \forall i\neq j, \widehat{H}_{ij}\neq0,
\\
 \widehat{Z}_{ij}  &\in & [-1,1] \qquad \forall i\neq j, \widehat{H}_{ij}=0,
\\
\widehat{H} &=&  \mathcal{P}_{\mathcal{F}^{k}} \biggl(\frac{1}{\tau} ( S-\rho
\widehat{Z} ) \biggr).
\end{eqnarray}

Let $\tilde H$, $\tilde Z$ be the support constrained FPS solution in
\eqref{eq:sparse-solution} and $\widehat{Z}$
be the dual variable constructed in \eqref{eq:dual-construct-1} to \eqref{eq:dual-construct-3}.
We first show that $(\tilde H,\widehat{Z})$ is also optimal for the
elastic net version of FPS when
$\tau$ is small enough.

From the existence proof above and Lemma~\ref{lem:primal-dual-optimality},
we know that (i) $S-\rho\widehat{Z}= \tilde\Sigma$,
(ii) the $k$-dimensional principal subspace of
$\tilde\Sigma$ is spanned by ${\hat U_J\choose 0}$ and (iii)~$\lambda_k(\tilde\Sigma)-\lambda_{k+1}(\tilde
\Sigma)>0$.

By the construction of $\tilde H$, part 3 of Lemma~\ref{lem:fantop-projection} implies
that when
%
%e32 #&#
\begin{equation}
\label{eq:tau-range}
0<\tau\le\lambda_k(\tilde\Sigma)-\lambda_{k+1}(
\tilde\Sigma)
\end{equation}
we have
\[
\tilde H=\mathcal{P}_{\mathcal{F}^{k}} \biggl(\frac{1}{\tau} (S-\rho\widehat{Z}
) \biggr).
\]

As a consequence, $(\tilde H,\widehat{Z})$ is also an optimal
primal--dual pair for the elastic net FPS
problem \eqref{eq:fps-en} when $\tau$ is in the range specified in
\eqref{eq:tau-range}.

Now we prove uniqueness of $\tilde H$ as a solution to the FPS problem
\eqref{eq:fps}. Assume that
there is another solution $\widehat{H}'\in\mathcal{F}^{k}$ such that
\[
\langle S, \tilde H \rangle-\rho\lVert\tilde H\rVert_{1,1} = \bigl
\langle S, \widehat{H}' \bigr\rangle-\rho\bigl\lVert
\widehat{H}'\bigr\rVert_{1,1}.
\]
But $\tilde H$ is the unique solution to the elastic net FPS for $\tau
>0$ small enough, we must have
$\lVert\widehat{H}'\rVert_F^2>\lVert\tilde H\rVert_F^2$, and hence
\[
k\ge\bigl\lVert\widehat{H}'\bigr\rVert_F^2>
\lVert\tilde H\rVert_F^2=k,
\]
which is a contradiction (the first inequality follows from that
$\widehat{H}'\in\mathcal{F}^{k}$).

%s6.1.3 #&#
\subsubsection{False negative control}
The false negative control under condition \eqref{eq:signal-strength-1} is obvious in view
of the Frobenius norm error bound \eqref{eq:frobenius-error}.

Now we prove false negative control under the entry-wise condition
\eqref{eq:entry-wise-min}.
According to Theorem~\ref{thm:deterministic}, we know that
$\widehat{H}$ is supported on $J$, and $\widehat{H}_{JJ}$ corresponds
to the
projector
of the $k$-dimensional principal subspace of $\tilde\Sigma
_{JJ}:=
S_{JJ}-\rho\widehat Z_{JJ}$
where $\widehat Z$ is the optimal dual variable.

Then it is
sufficient to
show that the leading eigenvector of $\tilde\Sigma_{JJ}$
does not have zero entries. % , which is equivalent to show that the
%leading eigenvector of
% $B\tilde\PopMatrix_{JJ} B$ does not have zero entries.
Note that $\lVert\tilde\Sigma_{JJ}-\Sigma_{JJ}\rVert_{\infty
,\infty}
\le2\rho$ and the second part of assumption~\eqref{eq:entry-wise-min}
implies that $\operatorname{sign}(\tilde\Sigma_{JJ})=\operatorname
{sign}(\Sigma_{JJ})$.

By the first part of assumption \eqref{eq:entry-wise-min},
we have $\operatorname{sign}(\tilde\Sigma_{JJ})=\operatorname
{sign}(\Sigma_{JJ})=b b^T$, where $b\in\{-1,1\}^s$.
Let $B$ be the $s\times s$ diagonal matrix such that $\operatorname
{diag}(B)=b$.
The matrix $B \tilde\Sigma_{JJ} B$ has all positive entries, and hence
by the Perron--Frobenius theorem, it has a unique leading eigenvector\vspace*{1pt}
$v_1$ whose entries are
all positive. As a result, the leading eigenvector of $\tilde\Sigma
_{JJ}$ is $Bv_1$, which does not
have zero entries.

%s6.1.4 #&#
\subsubsection{Auxiliary lemmas}

%le2 #&#
\begin{lemma}
\label{lem:sparse-solution}
Under the assumptions in Theorem~\ref{thm:deterministic}, let $\tilde H$ be
the solution
to the further constrained problem \eqref{eq:sparse-solution}.
Then $\tilde H$ is rank $k$ and unique.
Furthermore, there
exist $s\times k$ orthonormal matrices $U_J$, $\hat U_J$ such that:
\begin{enumerate}
\item ${U_J\choose 0}$ and ${\hat U_J\choose 0}$ span the $k$-dimensional principal subspaces of $\Sigma$ and
$S-\rho\tilde Z$,
respectively.
\item
There exists a $s\times s$ orthonormal matrix $Q$ such that
\begin{eqnarray*}
\hat U_J & = & Q U_J,
\\
\lVert Q-I \rVert_F & \le&\frac{8\rho s}{\lambda_k-\lambda
_{k+1}}.
\end{eqnarray*}
\end{enumerate}
\end{lemma}

\begin{pf}
Consider $\tilde\Sigma_{JJ}:= S_{JJ}-\rho\tilde Z_{JJ}$.
We know that $\tilde H_{JJ}$ maximizes
$\langle\tilde\Sigma_{JJ}, H \rangle$ over all $H\in\mathcal{F}^{k}_{s}$
(the trace-$k$ Fantope of $\mathbb{R}^{s \times s}$). We argue that
$\tilde H_{JJ}$ is unique and has rank $k$.
By condition \eqref{eq:deterministic-condition-1} and $\tilde Z\in
\mathbb B_p$,
we have
$\lVert\tilde\Sigma_{JJ}-\Sigma_{JJ}\rVert_{\infty,\infty}\le
2\rho$,
and hence $\lVert\tilde\Sigma_{JJ}-\Sigma_{JJ}\rVert_F\le2\rho s$.
On the other hand, by the \eqref{SPS} condition it is straightforward
to verify that
the leading $k$ eigenvectors of $\Sigma_{JJ}$ are those of $\Sigma$
confined on $J$, and hence $\lambda_{\ell}(\Sigma_{JJ})=\lambda
_\ell$ for
all $1\le\ell\le k$.
Furthermore, since $\Sigma_{JJ}$ is a principal submatrix of $\Sigma
$, we have
$\lambda_{k+1}(\Sigma_{JJ})\le\lambda_{k+1}$.
Therefore,
\[
\lambda_k(\tilde\Sigma_{JJ})-\lambda_{k+1}(\tilde
\Sigma_{JJ}) \ge\lambda_k(\Sigma_{JJ})-
\lambda_{k+1}(\Sigma_{JJ})-4\rho s\ge \lambda_k-
\lambda_{k+1}-4\rho s>0
\]
by condition
\eqref{eq:deterministic-condition-2}.
The first claim follows from part 1 of Lemma~\ref{lem:fantop-projection}.

The second claim is trivial when $s=k$. Now we focus on the case $s>k$.
By the \eqref{SPS} condition we know that the unique $k$-dimensional
principal subspace of $\Sigma$
is spanned by ${ U_J\choose 0}$ where $U_J$ is a $s\times k$ orthonormal
matrix.
Then $U_J$ spans the $k$-dimensional principal subspace of $\Sigma_{JJ}$.

Using the fact that
$\lambda_{k}(\Sigma_{JJ})-\lambda_{k+1}(\Sigma_{JJ})\ge
\lambda_k-\lambda_{k+1}$,
and applying Proposition 2.2 in \citet{Vu:2013}, we can choose the
right rotations for the columns
of $\hat U_J$ and $U_J$ so
that
\[
\lVert\hat U_J - U_J\rVert_F \le\bigl
\lVert\hat U_J\hat U_J^T - U_J
U_J^T\bigr\rVert_F.
\]
Using Lemma 4.2 of \citet{Vu:2013} and Cauchy--Schwarz, we have
\[
\bigl\lVert\hat U_J\hat U_J^T-U_JU_J^T
\bigr\rVert_F\le\frac{2}{\lambda
_{k}-\lambda_{k+1}}\lVert\tilde\Sigma_{JJ} -
\Sigma_{JJ}\rVert _F\le \frac{4\rho s}{\lambda_{k}-\lambda_{k+1}}.
\]
The above two inequalities jointly imply that
\[
\lVert\hat U_J-U_J\rVert_F\le
\frac{4\rho s}{\lambda_k-\lambda
_{k+1}}.
\]

Now let $\hat V=(\hat U_J, \hat U_J^c)$ be an $s\times s$ orthonormal
matrix, and similarly $V=(U_J, U_J^c)$. One can show that, using the
same argument as above, $\hat U_J^c$
and $U_J^c$ can be chosen such that
\[
\bigl\lVert\hat U_J^c-U_J^c
\bigr\rVert_F\le\frac{4\rho s}{\lambda_k-\lambda
_{k+1}}.
\]
Let $Q=\hat V V^T$, then $Q U_J=\hat U_J$ and
\[
\lVert I-Q\rVert_F = \bigl\lVert(V-\hat
V)V^T\bigr\rVert_F=\lVert\hat V-V\rVert_F
\le\frac{8\rho
s}{\lambda_k-\lambda_{k+1}}.
\]
\upqed\end{pf}

%le3 #&#
\begin{lemma}
\label{lem:primal-dual-optimality}
Under the assumptions of Theorem~\ref{thm:deterministic}, let $(\tilde
H,\tilde Z)$
be the optimal primal--dual pair of the additionally constrained FPS
problem \eqref{eq:sparse-solution}.
Let
${\hat U_J\choose 0}$, ${U_J \choose 0}$, and $Q$ be defined as in Lemma~\ref{lem:sparse-solution}. Let
$\tilde\Sigma$ be defined
as in \eqref{eq:sigma-tilde}. Then
\[
\lambda_{k}(\tilde\Sigma)-\lambda_{k+1}(\tilde\Sigma)>0
\]
and
$\tilde H$ is the unique projector of the  $k$-dimensional principal
subspace of $\tilde\Sigma$.
\end{lemma}

\begin{pf}
We start from a decomposition of $\tilde\Sigma$ as follows:
%
%e33 #&#
\begin{eqnarray}
\tilde\Sigma & = &
\pmatrix{S_{JJ}-
\rho\tilde Z_{JJ} & Q \Sigma_{JJ^c}\vspace*{3pt}\cr
\Sigma_{J^c J}Q^T & \tilde W+\Sigma_{J^c J^c}}
\nonumber
\\
& =&\pmatrix{
S_{JJ}-\rho\tilde
Z_{JJ} - Q\Sigma_{JJ} Q^T + Q
\Sigma_{JJ} Q^T & Q\Sigma_{JJ^c}
\vspace*{3pt}\cr
\Sigma_{J^cJ}Q^T & \tilde W+\Sigma_{J^cJ^c}
}
\nonumber
\\[-8pt]
\label{eq:magic-matrix}
\\[-8pt]
\nonumber
&=& \pmatrix{S_{JJ}-\rho\tilde
Z_{JJ} - Q\Sigma_{JJ} Q^T & 0
\cr
0&\tilde W } + \pmatrix{
Q\Sigma_{JJ} Q^T & Q\Sigma_{JJ^c}
\vspace*{3pt}\cr
\Sigma_{J^cJ}Q^T & \Sigma_{J^cJ^c}}
\nonumber
\\
&=&\mbox{``noise''}+\mbox{``signal.''}\nonumber
\end{eqnarray}

It can be directly verified that
${\hat U_J \choose 0}$ spans the $d$-principal subspace of
%
%e34 #&#
\begin{equation}\label{eq:eigen-signal}
\pmatrix{
Q\Sigma_{JJ}
Q^T & Q\Sigma_{JJ^c}
\vspace*{3pt}\cr
\Sigma_{J^c J}Q^T & \Sigma_{J^c J^c}}
 =
\pmatrix{
Q & 0
\cr
0 & I}
  \times \Sigma \times
\pmatrix{
Q^T & 0
\cr
0 & I}.
\end{equation}
Moreover, \eqref{eq:eigen-signal} implies that the eigenvalues of the
``signal'' part in the decomposition \eqref{eq:magic-matrix} are the
same as
those of $\Sigma$.

To sum up, we have so far shown that
${\hat U_J\choose 0}$ spans the $k$-dimensional principal subspace of the signal
part, with
spectral gap $\lambda_k-\lambda_{k+1}$.

Next, we need to show that the $k$-dimensional principal subspace
remains unchanged after adding the ``noise'' part.

First, the block-diagonal structure of the ``noise'' matrix in \eqref{eq:magic-matrix} ensures
that ${\hat U_J\choose 0 }$ spans one of its $k$-dimensional spectral subspace (a $k$-dimensional spectral subspace of
a $p\times p$ symmetric matrix $A$ means that if $v$ is in this
subspace, then
$Av$ is also in this subspace).

Second, we show that twice the operator norm of the ``noise'' part is smaller
than the gap between $k$th and $(k+1)$th eigenvalues of $\tilde\Sigma
$, which is $\lambda_k-\lambda_{k+1}$.
In fact,
the operator norm of the noise part
does not exceed
\begin{eqnarray*}
\lVert S_{JJ}-\rho\tilde Z_{JJ}-\Sigma_{JJ}
\rVert+\bigl\lVert\Sigma _{JJ} - Q\Sigma_{JJ} Q^T
\bigr\rVert &\le & 2\rho s+2\lVert\Sigma_{JJ}\rVert \times\lVert
Q-I\rVert
\\
&\le &  2\rho s + 2\lambda_1\times8\rho s/(\lambda_k-
\lambda_{k+1}),
\end{eqnarray*}
where the bound on $\lVert Q-I\rVert$ comes from Lemma~\ref{lem:sparse-solution}.
We also have
$\lVert\tilde W\rVert\le\rho$, which is contained within the above bound.

Therefore, by standard perturbation theory such as Weyl's inequality,
the subspace spanned by
${\hat U_J\choose 0 }$ is the $k$-dimensional principal subspace of $\tilde\Sigma$
as long as
%
%e35 #&#
\begin{equation}
\label{eq:eigen-gap-final} 4\rho s+16\sqrt{2}\lambda_1\rho s/(
\lambda_k-\lambda_{k+1})\le \lambda_k-
\lambda_{k+1} ,
\end{equation}
which means that twice the noise operator norm does not exceed the
spectral gap in the signal part.

When the inequality in \eqref{eq:eigen-gap-final} is strict, as stated
in condition \eqref{eq:deterministic-condition-2}, we know that the
$k$-dimensional principal subspace of
$\tilde\Sigma$ is unique.
\end{pf}

%s6.2 #&#
\subsection{Other proofs}\label{sec:additional-proof}
\mbox{}
\begin{pf*}{Proof of Lemma~\ref{lem:fantop-projection}}
(1) See \citet{OW:1992}.
(2) is Lemma~4.1 of \citet{fps-nips}.
(3) We have
\[
\langle A, H \rangle- \frac{\tau}{2} \lVert H\rVert_F^2
= - \frac{\tau}{2} \bigl\lVert H- \tau^{-1}A\bigr
\rVert_F^2 + \frac{1}{2\tau} \lVert A
\rVert_F^2.
\]
This is maximized over $H\in\mathcal{F}^{k}$ by
$H= \mathcal{P}_{\mathcal{F}^{k}}(\tau^{-1}A)$.
Note that by assumption $\gamma_{k} / \tau\geq1$
and $\gamma_{k+1} < \gamma_{k}$.
Then the claim follows by applying (1) and (2).
\end{pf*}

\begin{pf*}{Proof of Theorem~\ref{thm:agnostic-fps}}
Let $H_R$ be any solution of
\[
\max_{H\in\mathcal{F}^{k}, \lVert H\rVert_{1,1} \leq R} \langle\Sigma, H \rangle.
\]
Then $0 \leq\langle-\Sigma, \widehat{H}- H_R \rangle$,
and \eqref{eq:fps-constraint} implies
$0 \leq\langle S, \widehat{H}- H_R \rangle$.
Combining these two inequalities with the H\"older and
triangle inequalities yields
\[
0 \leq \langle\Sigma, H_R \rangle - \langle\Sigma, \widehat{H}
\rangle \leq \langle S- \Sigma, \widehat{H}- H_R \rangle \leq 2 R
\lVert S- \Sigma\rVert_{\infty,\infty}.
\]
Finally, invoke (2) to complete the proof.
\end{pf*}
\end{appendix}

\section*{Acknowledgments}
We thank the Editors and referees for their helpful comments.

% imsref loaded by daiva.urboniene, 2014-11-28 09:15:59


%suskaldyti doi

\printaddresses
\end{document}